\documentclass[12pt]{amsart}
\usepackage{amsmath}
\usepackage{amsfonts}
\usepackage{amsthm, upref}
\usepackage{graphicx}
\usepackage[usenames, dvipsnames]{color}
\usepackage{enumerate}

\usepackage{bm}

\bmdefine\betab{\mathbf{\beta}}
\bmdefine\sigmab{\mathbf{\sigma}}

    \numberwithin{equation}{section}

\newcommand{\comment}[1]{}
\newcommand{\eq}{\begin{equation}}
\newcommand{\en}{\end{equation}}
\newcommand{\rr}{\mathbb{R}}
\newcommand{\pp}{\mathbb{P}}

\newcommand{\ev}{\mathbb E}  

\newcommand{\gen}{\mathcal{L}}

\newcommand{\ep}{\hfill $\Box$}

\begin{document}

\theoremstyle{plain}
\newtheorem{thm}{Theorem}
\newtheorem{lemma}[thm]{Lemma}
\newtheorem{prop}[thm]{Proposition}
\newtheorem{cor}[thm]{Corollary}

\theoremstyle{definition}
\newtheorem{defn}{Definition}
\newtheorem{asmp}{Assumption}
\newtheorem{notn}{Notation}
\newtheorem{prb}{Problem}

\theoremstyle{remark}
\newtheorem{rmk}{Remark}
\newtheorem{exm}{Example}
\newtheorem{clm}{Claim}

\title[One-sided Tanaka equation]{On the one-sided Tanaka equation with   drift}

\author{Ioannis Karatzas, Albert N. Shiryaev and Mykhaylo Shkolnikov}
\address{INTECH Investment Management  \\ ONE PALMER SQUARE\\ PRINCETON, NJ 08542 and COLUMBIA UNIVERSITY \\ DEPARTMENT OF MATHEMATICS\\ NEW YORK, NY 10027}
\email{ik@enhanced.com, ik@math.columbia.edu}
\address{STEKLOV MATHEMATICAL INSTITUTE  \\ GUBKIN STREET 8 \\ 119991 MOSCOW, RUSSIA}
\email{albertsh@mi.ras.ru} 
\address{INTECH Investment Management  \\ ONE PALMER SQUARE\\ PRINCETON, NJ 08542 and STANFORD UNIVERSITY \\ DEPARTMENT OF MATHEMATICS\\ STANFORD, CA 94305}
\email{mshkolni@math.stanford.edu}

\keywords{Stochastic differential equation, weak existence, weak uniqueness, strong existence, strong uniqueness, Tanaka equation, skew Brownian motion, sticky Brownian motion, comparison theorems for diffusions}

    \subjclass[2000]{60H10, 60J60, 60J65}


\date{\today}

\begin{abstract} 
We study   questions of existence and uniqueness of weak and strong solutions  for a one-sided Tanaka equation with constant drift $\lambda\,$. We observe a dichotomy in terms of the values of the drift parameter:   for $\lambda \leq 0\,$, there exists a strong solution which is pathwise unique, thus also unique in distribution;  whereas 
for $\lambda>0\,$, the equation has a unique in distribution weak solution, but no strong solution (and not even a weak solution that spends zero time at the origin).  We also show that strength and pathwise uniqueness are   restored to the equation via suitable ``Brownian perturbations". 
\end{abstract}

\maketitle

\section{Introduction}

This paper   studies the   one-dimensional stochastic differential equation 
\eq
\label{sde}
\mathrm{d}X(t)=\lambda\;\mathrm{d}t + \mathbf{ 1}_{\{X(t)>0\}}\;\mathrm{d}W(t)\,,\qquad 0 \le t < \infty\,,
\en
where $W$ is  standard Brownian motion and $\lambda$  a real constant. The diffusion function $\, \sigmab (x) = \mathbf{ 1}_{ (0, \infty)} (x) \,$ in this equation is both discontinuous and degenerate, so questions of existence and uniqueness of   solutions are not covered by the classical theories of It\^o, Stroock \& Varadhan or Yamada \& Watanabe (e.g., Chapter 5 of \cite{KS}). 

\smallskip
When $\lambda =0\,$, the equation \eqref{sde} can be   viewed as a one-sided version of the {\it Tanaka  equation}
\eq
\label{Tanaka}
\mathrm{d}X(t)\,=\,\text{sgn}(X(t))\;\mathrm{d}W(t),\quad 0 \le t < \infty \,,
\en   
where  the signum function is defined as $\text{sgn} (x) = 1$ for $x >0$ and $\text{sgn} (x) = -1$ for $x \le0\,$.
 It was shown by Zvonkin \cite{Z} (e.g., Example 3.5, Chapter 5 of \cite{KS}) that the   equation (\ref{Tanaka}), for which weak existence and weak uniqueness (i.e., uniqueness in distribution) both hold, does not admit a strong solution; and that   strong (that is, pathwise) uniqueness fails for (\ref{Tanaka}).  
 
 \smallskip
 The equation (\ref{Tanaka}) is a special case of the {\it Barlow equation} 
 \eq
\label{Bar}
\mathrm{d}X(t)\,=\,  \alpha\;\mathbf{ 1}_{\{X(t)>0\}}\;\mathrm{d}W(t) - \beta\;\mathbf{ 1}_{\{X(t)\leq0\}}\;\mathrm{d}W(t)\,,\quad 0 \le t < \infty  
\en
with real  constants $\,\alpha >0\,$, $\beta >0\,$; as was shown by Barlow \cite{B}, for this equation weak existence and weak uniqueness hold  but strong uniqueness fails.

 \smallskip
At the same time, one can view a solution to the equation \eqref{sde} as a degenerate version of the {\it skew-Brownian motion} studied by Walsh \cite{W} and Harrison \& Shepp  
\cite{HS}, with the addition of a constant drift; see the recent paper \cite{ABTWW} and the survey \cite{L}, as well as the  references in these works. The skew-Brownian motion with constant drift is a solution of the equation
\eq
\label{skew}
\mathrm{d}X(t)\,=\,\lambda\;\mathrm{d}t + \alpha\;\mathbf{ 1}_{\{X(t)>0\}}\;\mathrm{d}W(t) + \beta\;\mathbf{ 1}_{\{X(t)\leq0\}}\;\mathrm{d}W(t)\,,\quad 0 \le t < \infty 
\en
for some given  real constants $\,\lambda  \,$, $\,\alpha >0\,$, $\beta >0\,$. It follows from the results of Nakao \cite{N} that this equation has a pathwise unique, strong solution. 

Formally letting $\, \beta \downarrow 0\,$ in (\ref{Bar}) and ``arguing by analogy", one might conjecture (as we did initially) that for the equation (\ref{sde}) with $\, \lambda =0\,$ strong existence and strong uniqueness fail. Similarly, letting $\, \beta \downarrow 0\,$ in (\ref{skew}) with $\, \lambda >0\,$, one might conjecture   that  the equation (\ref{sde}) with $\, \lambda >0\,$ has   a pathwise unique,  strong solution. Both conjectures would be wrong, a fact that illustrates the pitfalls of this kind of spurious reasoning. 

As it turns out, and as we   show below, for $\lambda \leq 0\,$  there exists a strong solution which is pathwise unique, thus also unique in distribution (Theorems \ref{Thm1} and \ref{Thm2}). With $\lambda>0\,$  the equation has a  weak solution which is unique in distribution, but has no strong solution (Theorem \ref{Thm3}); whereas not even a weak solution exists under the ``non-stickiness condition" $\,  Leb \, \{ t \ge 0 : X(t) = 0\}=0\,$ of \cite{MS}, where $\textit{Leb}\,$ stands for the Lebesgue measure on $[0,\infty)$ (Theorem \ref{Thm4}). The results for the equation (\ref{sde}) with $\, \lambda >0\,$ extend also to the more general equation 
$$
\mathrm{d}X(t)\,=\, \kappa \, \mathbf{ 1}_{\{X(t)>0\}}\;\mathrm{d}t + \lambda\, \mathbf{ 1}_{\{X(t)\leq0\}}\;\mathrm{d}t 
+ \mathbf{ 1}_{\{X(t)>0\}}\;\mathrm{d}W(t)\,,\qquad 0 \le t < \infty \,,
$$
for arbitrary $\, \kappa \in \mathbb{R}\,$ and $\, \lambda >0\,$. 

When $\lambda>0\,$, we show that suitable ``Brownian perturbations"  can restore to  the equation (\ref{sde}) a pathwise unique, strong solution  (Theorem \ref{Thm5}). 

 \medskip
 
\section{The case $\lambda<0$}

In the case $\lambda<0$ we show that the equation \eqref{sde} possesses a strong solution. Moreover, we prove pathwise uniqueness and uniqueness in distribution for the stochastic   equation \eqref{sde}.

\begin{thm}
 \label{Thm1}
Let $\lambda<0$. Then on each filtered probability space $(\Omega,\mathcal{F},(\mathcal{F}_t)_{t\geq0},\pp)$, which is rich enough to support a one-dimensional standard Brownian motion $W$ and a real-valued random variable $\zeta$, there exists a solution $X$ to the equation \eqref{sde}, which satisfies $X(0)=\zeta$ and is adapted to the filtration $\big(\mathcal {F}^{\,(\zeta, W)}_t\big)_{t\geq0}$ generated by $( \zeta, W)\,$. 

Moreover, $X$ is the unique process with these properties, and the distribution of every weak solution to   \eqref{sde} with the same initial distribution must coincide with the distribution of $X$.  
\end{thm}

 \smallskip

\noindent{\it Proof.} {\sl Step (A):} Considering a filtered probability space as posited in the statement of the theorem, we first define the stopping time
\eq
\label{tau}
\tau:=\inf\{t\geq0:\;\zeta+\lambda t+W(t)\leq 0\}
\en
and claim that the process $X(t)=\zeta + \lambda t + W( t \wedge \tau)$, $t\geq0$ is a strong solution of the equation \eqref{sde}. Indeed, $X(0)=\zeta$, the process $X$ is adapted to the filtration $\big(\mathcal {F}^{\,(\zeta, W)}_t\big)_{t\geq0}$,  and   for all  $\,0 \le t < \infty\,$ we have 
\begin{eqnarray*}
X(t)-X(0)=\lambda\, t + \int_0^t \mathbf{ 1}_{\{\tau>s\}}\;\mathrm{d}W(s)&=&\lambda t + \int_0^t \mathbf{ 1}_{\{\zeta+\lambda s + W( s \wedge \tau)>0\}} \;\mathrm{d}W(s)\\
&=&\lambda t + \int_0^t \mathbf{ 1}_{\{X(s)>0\}}\;\mathrm{d}W(s) \,.
\end{eqnarray*} 
 
 \medskip
\noindent {\sl Step (B):} Now, we claim that if $(Y,B,\xi)$ is a weak solution of the stochastic integral equation 
\[
 Y(t)\,=\, \xi + \lambda\,t + \int_0^t \mathbf{ 1}_{\{Y(s)>0\}}\;\mathrm{d}B(s)\,,\qquad 0 \le t < \infty
 \]
on an appropriate filtered probability space, then  $Y(t)=\xi + \lambda t + B( t \wedge \sigma)$ must hold for $\, 0 \le t < \infty\,$, where we have set 
\eq
\label{sigma}
\sigma\,:=\,\inf\{t\geq0:\;\, Y(t)\leq 0\}.
\en
This will immediately imply weak uniqueness and pathwise uniqueness for the stochastic differential equation \eqref{sde}.

To prove this claim, we fix a weak solution $(Y,B,\xi)$ with the described properties. Moreover, for every $\varepsilon>0$ we introduce the stopping times
\begin{eqnarray}
&&\sigma_{-\varepsilon}:=\inf\{t\geq0:\;Y(t)\leq -\varepsilon\},
\label{sigmae}\\
&&\varrho_{-\varepsilon }:=\inf\{t\geq\sigma_{-\varepsilon}:\;Y(t)\geq0\}. 
\end{eqnarray}
Suppose that for some $\varepsilon>0$ we had $\,\varrho_{-\varepsilon }<\infty\,$ on a set of positive probability; then on this same set  
\begin{eqnarray*}
0<\varepsilon=Y(\varrho_{-\varepsilon })-Y(\sigma_{-\epsilon})
=\lambda(\varrho_{-\varepsilon}-\sigma_{-\epsilon})+\int_{\sigma_{-\varepsilon}}^{\varrho_{-\varepsilon}} \mathbf{ 1}_{\{Y(s)>0\}}\;\mathrm{d}B(s)\\
=\lambda(\varrho_{-\varepsilon}-\sigma_{-\varepsilon})<0 
\end{eqnarray*}
would hold as well, which is clearly   absurd. This shows that  $\,\varrho_{-\varepsilon}=\infty\,$  is valid with probability one for all $\varepsilon>0$. Using this fact and
\[
Y(  \sigma_{-\varepsilon} \vee t)-Y(\sigma_{-\varepsilon})\,
=\,\lambda \big( ( \sigma_{-\varepsilon} \vee t )-\sigma_{-\varepsilon}\big)+\int_{\sigma_{-\varepsilon}}^{\, \sigma_{-\varepsilon} \vee t } \mathbf{ 1}_{\{Y(s)>0\}}\;\mathrm{d}B(s)\,,
\]
we conclude that 
\eq
\label{max}
Y(\sigma_{-\varepsilon} \vee t)\,=\,Y(\sigma_{-\varepsilon})+\lambda \big( ( \sigma_{-\varepsilon} \vee t )-\sigma_{-\varepsilon}\big) 
\en
holds for every $\,0 \le t < \infty\,$ and $\,\varepsilon>0\,$, with probability one.

On the other hand, we have
\[
 Y( t \wedge \sigma)\,=\,\xi+\lambda\cdot (t \wedge \sigma)+\int_0^{\,t \wedge \sigma} \mathbf{ 1}_{\{Y(s)>0\}}\;\mathrm{d}B(s)~~~~~~~~~~~~~
\]
\eq
\label{min}
 \qquad \qquad \qquad =\,\xi+\lambda\cdot (t \wedge \sigma)+B( t \wedge \sigma)\,  , \qquad 0 \le t < \infty\,,
\en 
therefore also
\eq
\label{sigmatoo}
\sigma\,=\,\inf\{t\geq0:\;\xi+\lambda t+B(t)\leq 0\}\,.
\en
 
Now, we claim, it is enough to show that the identity
\eq\label{stopptimes}
\lim_{\varepsilon\downarrow0} \,\sigma_{-\varepsilon}\,=\,\sigma 
\en
holds with probability one, in the notation of (\ref{sigma}), (\ref{sigmae}); because then, using \eqref{stopptimes} in conjunction with the   observations (\ref{min}), (\ref{max}) and the continuity of $\,   Y(\cdot)\,$, we will be able   to conclude
\begin{eqnarray*}
    Y(t) &=& Y(t)\cdot \mathbf{ 1}_{\{t\leq \sigma \}} + Y(t)\cdot \mathbf{ 1}_{\{t>\sigma\}}\\
&=&(\xi+\lambda\cdot  (t \wedge \sigma)+B( t \wedge \sigma))\;\mathbf{ 1}_{\{t\leq\sigma\}}+(Y(\sigma)+\lambda(t-\sigma))\;\mathbf{ 1}_{\{t>\sigma\}}\\
&=&(\xi+\lambda\cdot (t \wedge \sigma)+B( t \wedge \sigma))\;\mathbf{ 1}_{\{t\leq \sigma\}}+(\xi+\lambda \sigma+B(\sigma)+\lambda(t-\sigma))\;\mathbf{ 1}_{\{t>\sigma\}}\\
&=&\xi+\lambda t+B( t \wedge \sigma) 
\end{eqnarray*}
for all $\, 0 \le t < \infty\,$, as posited. 

\bigskip
\noindent {\sl Step (C):} We start   by recalling the Dambis-Dubins-Schwarz Theorem (cf.$\,$Theorem 4.6 and Problem 4.7, in Chapter 3 of \cite{KS}), according to which there is a one-dimensional standard Brownian motion $\,\betab (\cdot)\,$ such that
\eq
 \label{DDS}
M(t) \,:=\, \int_0^t \mathbf{ 1}_{\{Y(s)>0\}}\;\mathrm{d}B(s)\,=\,\betab\big(\langle M\rangle(t)\big)\,,\quad 0 \le t < \infty
\en 
holds, where   
$$
\langle M\rangle ( t) \,=\, \int_0^t \mathbf{ 1}_{\{Y(s)>0\}}\;\mathrm{d}s \,=\, Leb\,\big(  \{ s \in [0,t] : Y(s)>0\} \big)\, , \qquad \, 0 \le t < \infty\,
$$ 
 is the quadratic variation process of the martingale $\,M$ defined in (\ref{DDS})  and $\textit{Leb}\,$ stands for the Lebesgue measure on $[0,\infty)$.

To prove \eqref{stopptimes}, it suffices to show   $\, \mathbb{P} (E_\delta) =0\,$   for every $\,\delta>0\,$, where we introduce the  event 
\eq\label{deltaclaim}
E_\delta \,:=\, \big\{ M(\sigma+t) -  M(\sigma) \geq-\lambda \, t\,,\,\,\;\forall ~t\in[0,\delta] \,\big\}  
\en
with the notation of (\ref{DDS}).  In order to prove this   assertion, we fix a number  $\,\delta >0\,$, recall the notation (\ref{sigma}), and  define the  random  set 
\eq
A_\delta\,:=\, \big\{s\in[\,\sigma,\sigma+\delta\,]:\;Y(s)>0 \big\}\,.
\en

On the event  $\,  \{\textit{Leb}\,(A_\delta )=0\}  \,$ we have $\,\langle M\rangle ( \sigma + \delta ) = \langle M\rangle ( \sigma)\,$, thus $\,   M  ( \sigma + t ) =   M  ( \sigma)\,$ for all $\, t \in [0, \delta]\,$,  and consequently  $\, \mathbb{P} \,\big( E_\delta \cap \{\textit{Leb}\,(A_\delta )=0\} \big)\, =\,0\,$.    Thus, it suffices to show that the probability of the intersection of the event in \eqref{deltaclaim}  with the event $\,\{\textit{Leb}\,(A_\delta )>0\}\,$, is also equal to zero:
\eq
 \label{EA}
 \mathbb{P} \,\big( E_\delta \cap \{\textit{Leb}\,(A_\delta )>0\} \big)\, =\,0\,.
 \en

In order to do this, we use the representation (\ref{DDS}) to see that the event in \eqref{deltaclaim} is contained in
\[
\widetilde{E}_\delta  :=  \bigg\{\inf_{ \langle M\rangle(\sigma)\,\le \,s\, \le \, \langle M\rangle(\sigma)+\textit{Leb}\,(A_\delta)}
\Big(\lambda \big( \langle M\rangle^{-1}(s)-\langle M\rangle^{-1}(\sigma)\big)+\betab(s)-\betab \big(\langle M\rangle(\sigma)\big)\Big)\geq 0\bigg\};
\]

\medskip
\noindent
 here we have remarked $\,\langle M\rangle(\sigma+ \delta) =  \langle M\rangle(\sigma)+\textit{Leb}\,(A_\delta)\,$ and set $$\langle M\rangle^{-1}(s)\,:=\inf\{t\geq0:\;\langle M\rangle(t)>s\}\,, \qquad s\geq0\,.$$ 
On the intersection of events $\,  \widetilde{E}_\delta \cap \{\textit{Leb}\,(A_\delta )>0\}  \,$, the one-dimensional, standard Brownian motion
$$
\widetilde{\betab} (u) \,:=\,\betab \big(\langle M\rangle(\sigma) +u\big) - \betab \big(\langle M\rangle(\sigma)\big)\,, \qquad 0 \le u < \infty
$$
has to stay nonnegative throughout the interval $\, [\,0\,,  \,Leb \,(A_\delta)\,]\,$ (recall here that $\, \lambda <0$); but this implies $\, \mathbb{P} \,\big( \widetilde{E}_\delta \cap \{\textit{Leb}\,(A_\delta )>0\} \big)=0\,$, because a one-dimensional, standard Brownian motion changes sign infinitely often on every non-trivial time interval starting at the origin, with probability one (e.g., Problem 7.18 in Chapter 2 of \cite{KS}). In conjunction with the inclusion $\, E_\delta \subseteq \widetilde{E}_\delta\,$, this observation leads to (\ref{EA}) and completes the proof. \ep  

\section{The case $\lambda=0$}
 
In the case $\lambda=0$ we prove the existence of a strong solution to the equation \eqref{sde}. Moreover, we show pathwise and weak uniqueness for the equation \eqref{sde}. Thus, the one-sided Tanaka equation   $\, \mathrm{d} X(t) = \mathbf{ 1}_{ \{ X(t) >0 \} }\, \mathrm{d} W(t)\,$ of \eqref{sde}  has  qualitative properties markedly different from those of the ``real" Tanaka equation \eqref{Tanaka}. 

We remark at this point that this can be shown along the lines of the proof in the case $\lambda<0$, but we  prefer to give here a shorter proof  which relies on the Engelbert-Schmidt \cite{ES} criterion for weak uniqueness   of one-dimensional stochastic differential equations without drift. 

\begin{thm}
 \label{Thm2}
Let $\,\lambda=0\,$. Then, on each filtered probability space $(\Omega,\mathcal{F},(\mathcal{F}_t)_{t\geq0},\pp)$  rich enough to support a one-dimensional standard Brownian motion $W$ and an independent real-valued random variable $\zeta$, there exists a solution $X$ to the equation \eqref{sde}, which satisfies $X(0)=\zeta$ and is adapted to the filtration $\big(\mathcal {F}^{\,(\zeta, W)}_t\big)_{t\geq0}$ generated by $(\zeta, W)\,$. 

Moreover, $X$ is the unique process with these properties, and the law of every weak solution to   \eqref{sde} with the same initial distribution must coincide with the law of $X$.  
\end{thm}

\noindent{\it Proof:} We start by defining the stopping time
\eq
\vartheta\,:=\,\inf\{t\geq0:\;\zeta+W(t)\leq0\} 
\en  
and setting $X(t)=\zeta+W(  t \wedge \vartheta)$, $\,t\geq0$. Arguing as in Step (A) in the proof of Theorem \ref{Thm1}, we show that $X$ is a strong solution of the stochastic differential equation \eqref{sde} with initial value $\zeta$. 

\smallskip
We claim now that weak uniqueness holds for the equation \eqref{sde} with $\lambda=0$. To this end, we employ the Engelbert \& Schmidt \cite{ES} theory in the form of Theorem 5.7 in Chapter 5 of \cite{KS}, and need to show the following identity between sets:
\eq
\label{Schmidt}
\big\{x\in\rr:\;  \sigmab (x)=0\big\}\,=\,\Big\{x\in\rr:\;\int_{x-\epsilon}^{x+\epsilon} \frac{\mathrm{d}y}{\,\sigmab^2(y)\,}=\infty\,,\;\;~\forall~~\epsilon>0\Big\}.
\en
Indeed, one checks fairly easily that for the diffusion function $\, \sigmab (x) = \mathbf{ 1}_{(0,\infty)} (x)\,$ of the equation (\ref{sde}), both sets in (\ref{Schmidt}) are equal to  $(-\infty,0]\,$, so that weak uniqueness holds for the equation \eqref{sde} with $\lambda=0$.

 \smallskip
It remains to show strong uniqueness. To this end, let $Y$ be another solution of the equation \eqref{sde} defined on the same probability space as $X$, adapted to the filtration $\big(\mathcal {F}^{\,(\zeta, W)}_t\big)_{t\geq0}$ and satisfying $Y(0)=\zeta=X(0)$. From the explicit formula for the process $X$ we deduce
\eq
\vartheta=\inf\{t\geq0:\,X(t)\leq0\}.
\en 
We now define a new stopping time $\vartheta'$ by
\eq
\vartheta'\,:=\,\inf\{t\geq0:\,Y(t)\leq0\}\,.
\en 
Due to weak uniqueness, we must have with probability one: $Y(t\vee\vartheta')=X(t\vee\vartheta)=X(\vartheta)$ for all $\,0 \le t < \infty\,$. Moreover, from the stochastic differential equation \eqref{sde} we conclude
\[
Y(t\wedge\vartheta')=\zeta+W(t\wedge\vartheta'),\qquad 0 \le t < \infty\,.
\]
Thus, with probability one we have: $X(t\wedge\vartheta\wedge\vartheta')=Y(t\wedge\vartheta\wedge\vartheta')$, $\,~\forall\,\, \,0 \le t < \infty\,$. Combining the latter two observations we see that, in order to prove strong (pathwise) uniqueness,  it suffices to show that $\,\vartheta=\vartheta'\,$ holds with probability one. 

To this end, we note
\[
\vartheta'=\inf\{t\geq0:\,Y(t\wedge\vartheta')\leq0\}=\inf\{t\geq0:\,\zeta+W(t\wedge\vartheta')\leq0\}.
\]
This last expression is equal  to $\inf\{t\geq0:\,\zeta+W(t)\leq0\}=\vartheta\,$  if $\, \vartheta\leq\vartheta'$, and to infinity  if $\vartheta>\vartheta'$. However, the second case occurs with zero probability, since weak uniqueness implies $\,\pp(\vartheta'<\infty)=\pp(\vartheta<\infty)=1\,$. \ep

 \medskip

\section{The case $\lambda>0$}  

In the case $\lambda>0$ we shall show first that the equation \eqref{sde} has a unique weak solution, but  not   a strong solution. We shall also show that this solution is ``sticky at the origin", in the sense that the so-called {\sl non-stickiness condition} $\,\int_0^\infty \mathbf{ 1}_{\{X(t)=0\}} \, \mathrm{d} t =0\,$  cannot possibly hold with probability one. 

\medskip
\begin{thm}
  \label{Thm3}
Let $\,\lambda>0\,$. Then the equation \eqref{sde} has a   weak solution which is unique in the sense of the probability distribution, but does not admit a strong solution.
\end{thm}   

\noindent{\it Proof.}   {\sl Step (i):} We start by {\it constructing a weak solution} to the equation \eqref{sde}. In order to do this, we use the results of \cite{chi} (consult  also  pages 193-205 in the book \cite{GS},  and the more recent articles \cite{WJ1}, \cite{WJ2}) to conclude that  the equation
\eq\label{sticky'}
\mathrm{d}X(t)=\lambda\cdot\mathbf{ 1}_{\{X(t)\leq0\}}\;\mathrm{d}t + \mathbf{ 1}_{\{X(t)>0\}}\;\mathrm{d}B(t)\,,\qquad 0 \le t < \infty
\en
has a weak solution on a suitable filtered probability space $(\Omega,\mathcal{F},(\mathcal{F}_t)_{t\geq0}, \mathbb{Q} )$ for all initial values $X(0)\in\rr$, where $B$ is a one-dimensional standard Brownian motion. Indeed, for $X(0)\geq0$ one can define $X$ to be the sticky Brownian motion started at $X(0)$ (see \cite{chi}); and for $X(0)<0$, one can set $X(t)=X(0)+\lambda \, t$, $\,0\leq t\leq |X(0)|/\lambda$, then let $\,X(t)$, $\,  |X(0)|/\lambda \le t < \infty \,$ be a sticky Brownian motion started at the origin. 

\smallskip
Now  we carry out a Cameron-Martin-Girsanov change of probability measure,   from the underlying    $\,\mathbb{Q}\,$ to a probability measure $ \,\pp\, $   under which the process $\,W(t):=B(t)-\lambda \, t\,$, $\,0 \leq t < \infty\,$ is a standard Brownian motion. (See Corollary 5.2 in Chapter 3 of \cite{KS}, or pages 325-330 in \cite{RY}, for the details; the two measures $\,\mathbb{Q}\,$ and $\,\pp\,$ are equivalent  when restricted to $\,   \mathcal{F}_T \,$, for each $\, T \in (0, \infty)$.)  Under this new measure $\,\pp \,$, the process $X$ will satisfy 
\[
\mathrm{d}X(t)=\lambda\;\mathrm{d}t + \mathbf{ 1}_{\{X(t)>0\}}\;\mathrm{d}W(t)\,,\qquad 0 \le t < \infty\,,
\]
and  thus  $(X,W)$ will be a weak solution to \eqref{sde} on $(\Omega,\mathcal{F},(\mathcal{F}_t)_{t\geq0}, \pp )$.

\medskip
\noindent {\sl Step (ii):} Next, we prove {\it weak uniqueness.} To this end, let $(X,W)$ be an arbitrary weak solution of the equation \eqref{sde} on some probability space $(\Omega,\mathcal{F},(\mathcal{F}_t)_{t\geq0}, \pp )$. If we carry out again a Cameron-Martin-Girsanov change of measure such that  $\,B(t) =W(t)+\lambda\, t\,$, $\,\,0 \leq t < \infty\,$ becomes a standard Brownian motion under the new measure $\, \mathbb{Q}\,$
, then under this new measure the pair of processes $(X, B)$ will constitute a weak solution of the equation \eqref{sticky'}. 

\smallskip
We show now that, if the initial condition $\, X(0)\,$ is nonnegative, then the state process $\,X\,$ of such a weak solution remains nonnegative at all times. To do this, we pick     a nonincreasing function $f:\,\rr\rightarrow[0,1]$ supported in $(-\infty,0)$, which is twice continuously differentiable and has bounded first and second derivatives. Fixing   $\,t\in [0, \infty)\,$ and combining It\^o's formula with Fubini's Theorem, we deduce 
\[
\ev^{\,\mathbb{Q}} \,[\,f(X(t))\,]-f(X(0))=\lambda\cdot\int_0^t \ev^{\,\mathbb{Q}}\,\big[\,f'(X(s))\cdot \mathbf{ 1}_{\{X(s)<0\}}\,\big]\;\mathrm{d}s\leq 0\,,
\]
where $\, \ev^{\,\mathbb{Q}}\,$ denotes integration with respect to the auxiliary probability measure $\,  \mathbb{Q} \,$. In particular, we see that $X(0)\geq0$ implies $\ev^{\,\mathbb{Q}}\,[\,f(X(t))\,]=0\,$. Since the indicator function of every nonempty open interval in $(-\infty,0)$ can be dominated by a function $f$ as described above, and since the paths of $X$ are continuous, we conclude that $$X(0)\geq0 \qquad  \mathrm{implies} \qquad X(t)\geq0 \quad \hbox{for all}\quad  0 \le t < \infty\,,$$ with $\,  \mathbb{Q}-$probability one. Because the measures $\,\mathbb{P}\,$ and $\, \mathbb{Q}\,$ are equivalent  when restricted to $\,   \mathcal{F}_T \,$  for each $\, T \in (0, \infty)$, we see that the above implication holds also with $\,  \mathbb{P}-$probability one.

On the other hand, if $X(0)<0$, then the equation \eqref{sde} shows $\,X(t)=X(0)+\lambda t\,$ for all  $\,0\, \leq t\leq |X(0)|/\lambda\,$. The same argument as before, but now on the time interval $[\,|X(0)|/\lambda,\infty)\,$, yields $\,X(t)\geq0\,$ for all $\, t \in [\,|X(0)|/\lambda,\infty)\,$, with $\,  \mathbb{Q}-$probability one (thus also with $\,  \mathbb{P}-$probability one). 

We conclude that, under the new measure $\,\mathbb{Q}\,$, the process $X$ satisfies the stochastic differential equation 
\eq\label{sticky}
\mathrm{d}X(t)=\lambda\cdot\mathbf{ 1}_{\{X(t)=0\}}\;\mathrm{d}t + \mathbf{ 1}_{\{X(t)>0\}}\;\mathrm{d}B(t)  
\en
  driven by the $\,  \mathbb{Q}-$Brownian motion $\,B\,$ and therefore, on the strength of the results in \cite{chi}, \cite{WJ1}, has the distribution of the ``sticky Brownian motion" for all $\,t\in [0, \infty)\,$ if $\,X(0)\geq0\,$, and for all $t\in [\, |X(0)|/\lambda, \infty)\,$   if $\,X(0)<0\,$. Moreover, the main result in \cite{chi} shows that the joint distribution of the pair $(X,B)$ under $\mathbb{Q}\,$ is uniquely determined. Thus, making a change of measure back from $\mathbb{Q}$ to $\mathbb{P}$, we conclude that the distribution of $X$ under $\mathbb{P}$ must coincide with the distribution of the weak solution constructed above. This proves weak uniqueness.

 \medskip
\noindent {\sl Step (iii):}
Finally, we show by contradiction that the equation \eqref{sde} {\it cannot have a strong solution.} To this end, we suppose that $X$ is a strong solution to \eqref{sde} on a probability space $(\Omega,\mathcal{F},(\mathcal{F}^W_t)_{t\geq0},\pp)$; that is, $X$ solves \eqref{sde} and is adapted to the filtration $(\mathcal{F}^W_t)_{t\geq0}$ generated by the Brownian motion $W$ driving   the equation \eqref{sde}. Then, the same argument as in the proof of weak uniqueness shows that there is a Cameron-Martin-Girsanov change of measure, such that: the process  $B(t):=W(t)+\lambda \, t\,$, $\,0 \leq t < \infty\,$ is a standard Brownian motion under the new measure $\, \mathbb{Q}\,$; whereas   $X$ solves under this new measure the equation \eqref{sticky} for $\, 0 \le t < \infty\,$  if $\,X(0)\geq0\,$, and for $\,t\geq |X(0)|/\lambda\,$  if $\,X(0)<0\,$. But the processes   $\,W$ and   $B$ generate  exactly the same filtrations, so we conclude that for $X(0)\geq0$ we have constructed a strong solution of the equation \eqref{sticky}. This is in clear contradiction to the results in \cite{chi} and    \cite{WJ1}; Theorem 1 in the paper \cite{WJ1} shows, in particular,    that the conditional distribution of the sticky Brownian motion $\, X(t)$, given the entire path of the ``driving" Brownian motion $B$ in  (\ref{sticky}), is given by
$$
\mathbb{Q} \,\big( X(t) \le x \, \big|\, B(u)\,, \, 0 \le u < \infty \big) \,=\, \exp \big( - 2\, \lambda \, \big (B(t) + S(t) -x \big) \big)\,, \qquad 
$$ 
$\mathbb{Q}-$a.s., for all $\, x \in [0, B(t) + S(t)]\,$, where $\, S(t) := \max_{0 \le u \le t}  ( - B(u) )\,$.  Hence, a strong solution to the equation \eqref{sde} cannot exist. \ep 

\medskip
  Next, we provide a direct argument showing that for $\lambda>0$ the equation \eqref{sde} does not admit a weak solution which spends zero time at the origin (the ``non-stickiness  condition" (\ref{0timeat0}) below, in a terminology borrowed from  \cite{MS}). Clearly, this can  also be deduced from the weak uniqueness in Theorem \ref{Thm3}, and from the fact that the weak solution constructed in the proof of that result spends at the origin a non-zero amount of time with positive probability.  The method of proof of Theorem \ref{Thm4}, however, seems to be novel; it   might prove useful in the context of other stochastic differential equations, for which an analogue of Theorem \ref{Thm3} is not readily available. 

\begin{thm}
 \label{Thm4}
Let $\kappa$ be an arbitrary real constant and $\lambda$ an arbitrary positive constant. Then the stochastic differential equation 
\eq
\label{4.3}
~~~~~~~
\mathrm{d}X(t)\,=\, \kappa \, \mathbf{ 1}_{\{X(t)>0\}}\;\mathrm{d}t + \lambda\, \mathbf{ 1}_{\{X(t)\leq0\}}\;\mathrm{d}t 
+ \mathbf{ 1}_{\{X(t)>0\}}\;\mathrm{d}W(t)\,,\qquad 0 \leq t < \infty  \en
has no weak solution which satisfies
\eq
\label{0timeat0}
\int_0^T \mathbf{ 1}_{\{X(t)=0\}} \, \mathrm{d} t \,=\,0\,, \qquad  \mathrm{a.s.}
\en 
for all $\,T\in[0, \infty)\,$. 

In particular, this is the case for $\,\kappa=\lambda >0\,$, which corresponds to the equation \eqref{sde}.
\end{thm}

\noindent{\it Proof.} {\sl Step (1):} We shall suppose that $(X,W)$ is a weak solution of the equation \eqref{4.3} defined on a suitable filtered probability space $(\Omega,\mathcal{F},(\mathcal{F}_t)_{t\geq0},\pp)$ and satisfying \eqref{0timeat0}, and will derive a contradiction. 

To this end, we first carry out a Cameron-Martin-Girsanov change of probability measure,   from the underlying    $\,\pp$ to a probability measure $\,\widehat{\pp}\,$   under which the process $$\,\widehat{W}(t)\,:=\,W(t)+(\kappa+\lambda)\, t\, , \qquad \,0 \leq t < \infty\,$$ is a standard Brownian motion (the two probability measures $\,\widehat{\pp}\,$ and $\,\pp\,$ are equivalent  when restricted to $\,   \mathcal{F}_T \,$, for each $\, T \in (0, \infty)$). Substituting this into \eqref{4.3} we see that, under $\,\widehat{\pp}\,$, the process $X$ satisfies the stochastic differential equation 
\eq
\label{bang-bang}
\mathrm{d}X(t)\,=\,-\lambda\cdot \mathrm{sgn} \big( X(t)\big)\;\mathrm{d}t + \mathbf{ 1}_{\{X(t)>0\}}\;\mathrm{d}\widehat{W}(t)\,,
\en
with the distribution of $X(0)$ being unchanged. 

Next, we let $\,\mu(t)$, $\,0 \leq t < \infty\,$ be the collection of one-dimensional marginal distributions  of the process $X$ under $\widehat{\pp}$, namely $\, \mu(t) = \widehat{\mathbb{P}} \circ (X(t))^{\,-1}\,$. We shall show in Step 2 that the family of measures $\,\frac{1}{\,T\,}\int_0^T \mu(t)\,\mathrm{d}t\,$, $\,0 < T < \infty\,$ is uniformly tight, and in Step 3 that every limit point of this family must be the zero measure. This will  establish the desired contradiction.

\medskip

\noindent {\sl Step (2):} The uniform tightness of the family $\,\mu(t)$, $t \ge 0\,$ (and, hence, also of the family of measures $\,\frac{1}{\,T\,}\int_0^T \mu(t)\,\mathrm{d}t\,$, $\,\,0 < T < \infty$) will follow from the ``mean stochastic comparison" results of  Hajek \cite{H}. 

\smallskip
 
We start by recalling the definition $$L^{\Theta} (t) \,:=\, \Theta^+ (t) - \Theta^+ (0) - \int_0^t \mathbf{1}_{\{\Theta(s)>0\}} \, \mathrm{d} \Theta(s) \,=\,\lim_{\varepsilon \downarrow 0}\, { 1 \over \, 4 \, \varepsilon\,} \int_0^t \mathbf{ 1}_{\{|\Theta(s)| <\varepsilon\}}\,\mathrm{d} \langle \Theta \rangle (s)$$  of    the local time   accumulated  at the origin by a generic  continuous semimartingale $\,\Theta\,$  during the time-interval $[0,t]$, where  $\, \langle \Theta \rangle  \,$ is the quadratic variation of the local martingale part of $\,\Theta\,$. We recall also   the fact that the local time $\, L^\Theta (\cdot)\,$ is flat off the set $\, \{ t \in [0, \infty): \Theta(t)=0\}\,$; to wit, for every $\, T \in (0,\infty)\,$, we have   
\eq
\label{LT}
 \int_0^T \mathbf{ 1}_{\{\Theta(t)\neq 0\}}\;\mathrm{d}L^\Theta(t) \,=\,0  \,, \qquad \mathrm{a.s.}
\en
  (cf. Theorem 7.1, equation (7.2) on page 218, Chapter 3 of \cite{KS}).

\smallskip
With this terminology and notation in place, we  apply first the generalized It\^o rule (see Theorem 7.1, equation (7.4) on page 218, Chapter 3 of \cite{KS}) to the function $\, f(x) = |x|\,$ and the semimartingale $\,X\,$ of (\ref{bang-bang}), and obtain
\eq
\label{|bang-bang|}
\mathrm{d}\,|X|(t)\,=\,-\lambda\;\mathrm{d}t + \mathbf{ 1}_{\{X(t)>0\}}\;\mathrm{d}\widehat{W}(t) +2\, \mathrm{d}L^{X}(t)\,.
\en

 Next, we apply  the generalized It\^o rule, once again to the function $\, f(x) = |x|\,$ but this time to the semimartingale $ |X| $ of (\ref{|bang-bang|}); in conjunction with (\ref{LT}), we obtain
 
 $$
\mathrm{d}\,|X|(t)\,=\, \mathrm{sgn}\big( |X(t)| \big)  \big( 
-\lambda\;\mathrm{d}t + \mathbf{ 1}_{\{X(t)>0\}}\;\mathrm{d}\widehat{W}(t) +2\, \mathrm{d}L^{X}(t) \big) +2\, \mathrm{d}L^{|X|}(t) \qquad 
$$
$$
\qquad \qquad =\, - \lambda \big( 1 - 2\cdot \mathbf{ 1}_{\{X(t)=0\}} \big) \;\mathrm{d}t + \mathbf{ 1}_{\{X(t)>0\}}\;\mathrm{d}\widehat{W}(t) -2 \, \mathrm{d}L^{X}(t) +2\, \mathrm{d}L^{|X|}(t) \,.
$$

\medskip
\noindent
Comparing this last expression with (\ref{|bang-bang|}) and invoking the condition \eqref{0timeat0}, we deduce that for every $\, T \in (0,\infty)\,$ the identity
$$
L^{|X|} (T) \,=\, 2\, L^{X} (T)
$$ 
holds almost surely under both $\,\widehat{\pp}\,$ and $\,\pp\,$. Thus, the equation (\ref{|bang-bang|}) takes the form
\eq
\label{||bang-bang||}
\mathrm{d}\,|X|(t)\,=\,-\lambda\;\mathrm{d}t + \mathbf{ 1}_{\{X(t)>0\}}\;\mathrm{d}\widehat{W}(t) +   \mathrm{d}L^{|X|}(t)\,.
\en

Next, we consider the strong solution of the equation
\eq
\mathrm{d}Z(t)\,=\,-\lambda\;\mathrm{d}t + \mathrm{d}\widehat{W}(t) + \mathrm{d}L^{Z}(t)
\en 
satisfying the initial condition $Z(0)=|X(0)|$, where $L^{Z}(t)$ is the local time accumulated by $Z$ at the origin during the time-interval $[0,t]$. The process $Z$ can be constructed by applying the Skorohod map to the paths of the process $\,|X(0)|-\lambda t+\widehat{W}(t)\,$, $\,0 \leq t < \infty\,$. 

\smallskip
Since $Z$ is  Brownian motion with negative drift and   reflection at the origin, we know that it is a Markov process  with   (a unique) invariant distribution that has  exponential density $\, 2 \,\lambda\, e^{- 2 \lambda \, \xi}\,$, $\,\xi >0\,$. In particular, the family of one-dimensional marginal distributions of $Z$ is uniformly tight, and applying Theorem 2 of \cite{H} to the processes $|X|$ and $Z$  we obtain that, for every given $\varepsilon>0$, there exists a real number $K_\varepsilon>0$ such that  
\[
\sup_{0 \leq t < \infty}\,\mu(t)\big(\rr\backslash[-K_\varepsilon,K_\varepsilon]\big)\leq 2\cdot\sup_{0 \leq t < \infty}\,\widehat{\pp}\big(Z(t)\geq K_\varepsilon \big)<\varepsilon\,.
\]

We conclude that the family $\,\mu(t)$, $\,0 \leq t < \infty\,$ is uniformly tight. In particular, we can find a sequence $\,0<T_1<T_2<\dots\,$ of numbers that increase to infinity, for which the weak limit
\eq
 \label{WL}
\nu\,:=\,\lim_{n\rightarrow\infty}\,\,\frac{1}{\,T_n\,}\int_0^{T_n} \mu(t)\;\mathrm{d}t
\en
is well-defined and a probability measure on $\,\mathcal{B}  (\rr)\,$.

\medskip

\noindent {\sl Step (3):} We shall show now that the weak limit $\,\nu$ in (\ref{WL}) can only be the zero measure, and this will lead to the desired contradiction. We denote by $C_0^\infty(\rr)$ the space of continuous and infinitely continuously differentiable functions $\,f: \rr \rightarrow \rr\,$ which vanish at infinity together with all their derivatives. Applying It\^o's formula under the measure $\,\widehat{\pp}\,$, we see in conjunction with (\ref{bang-bang}) that the family of probability measures $\, \mu(t) = \widehat{\mathbb{P}} \circ (X(t))^{\,-1}\,$, $0 \leq t < \infty\,$ on $\,   \mathcal{B}  (\rr)\,$ satisfies the Fokker-Planck equation
\eq
\label{FP}
\forall \, ~f\in C_0^\infty(\rr),\,\;T \in (0, \infty):\quad \big(\mu(T),f\big)\,=\,\big(\mu(0),f\big)+\int_0^T \big(\mu(t), \gen f \big)\,\mathrm{d}t\,.
\en
Here we  denote by $(\cdot\,,\cdot)$ the pairing between finite measures and bounded measurable functions on $\rr\,$, and   have defined
\eq
\gen\,:=\,-\lambda\cdot \mathrm{sgn}(x)\, \frac{\mathrm{d}}{\mathrm{d}x} +\frac{1}{\,2\,}\cdot\mathbf{ 1}_{(0,\infty)}\, \frac{\mathrm{d}^2}{\mathrm{d}x^2}\,.
\en

\medskip
Now, we fix a constant $\,K>0\,$, and pick a function $\,f\in C_0^\infty(\rr)\,$ and a constant $b>2\lambda$ with the following properties: $f(x)=e^{\,b \,x}$, whenever $x\leq K\,$; $\,\gen f\geq1$, whenever $x\in[-K,K]\,$; and $\,\gen f\geq -1\,$, whenever $\,x\geq K\,$. This can be achieved by taking $b$ to be large enough first, and by choosing $f$ on the interval $[K,\infty)$ appropriately thereafter. Plugging $f$ into the Fokker-Planck equation (\ref{FP}) with $T=T_n$, dividing both sides of the equation by $T_n$ and taking the limit as $n\rightarrow\infty$, we get   
\[
0\,=\,\lim_{n\rightarrow\infty} \,\,\frac{\,1\,}{\,T_n\,}\int_0^{T_n} \big(\mu(t),\gen f\big)\,\mathrm{d}t\,.  
\]
Moreover, using the inequality
\[
\gen f \,\geq \,\mathbf{ 1}_{(-K,K)}-\mathbf{1}_{[K,\infty)},
\]
and applying the Portmanteau Theorem, we end up with
\[
0 \, \geq \, \nu \big((-K,K)\big)-\nu \big([K,\infty)\big).
\]
Hence, by taking the limit as $K\rightarrow\infty$, we obtain $\nu((-\infty,\infty))=0$, which provides the desired contradiction. \ep

\begin{rmk}
The result of Theorem \ref{Thm4} for the equation \eqref{sde}, that is, when $\kappa=\lambda\,$ in (\ref{4.3}), can be also obtained by the following shorter  but less instructive argument. 

We shall suppose that on a suitable filtered probability space $(\Omega,\mathcal{F},(\mathcal{F}_t)_{t\geq0},\pp)$ there is defined a weak solution $(X,W)$ of the equation \eqref{sde}, which   satisfies the non-stickiness  condition  \eqref{0timeat0}; and will derive a contradiction. For simplicity, we shall assume that $\, X(0) \,$ is a nonnegative constant. As in step (ii) in the proof of Theorem \ref{Thm3}, we conclude that, with probability one, we must have $\, X(t) \ge 0\,$ for all $\,0 \le t < \infty\,$, with probability one. Consequently,
\eq
\begin{split}
\label{sie}
 0 \, \le \, X(t)\,=\,X(0) + \lambda\; t + \int_0^t \mathbf{ 1}_{\{X(s)>0\}}\;\mathrm{d}W(s)\\ 
=\,X(0) + \lambda\; t + \int_0^t \mathbf{ 1}_{\{X(s)\ge 0\}}\;\mathrm{d}W(s)\,=\,X(0) + \lambda\; t +W(t)\,,\quad 0 \le t < \infty
\end{split}
\en
must also hold with probability one. The second equality in (\ref{sie}) is a consequence of $\, \int_0^\cdot \mathbf{ 1}_{\{X(s)=0\}} \, \mathrm{d}W(s) \equiv 0\,$, which is in turn a consequence of (\ref{0timeat0}); whereas the inequality and the third equality are  consequences of the nonnegativity of $X$. 

The inequality between the left- and right-most members in (\ref{sie}) implies that the probability of the event
$$
\big\{ \, \lambda \, t + W(t) < -X(0)\,, \quad \hbox{for some}~~\, ~~t \in [0, \infty)\, \big\}
$$
 is zero. We know, however (e.g., \cite{KS}, Exercise 5.9 on page 197), that the probability of this event is actually $\, e^{\, - 2 \, \lambda \, X(0)} >0\,$, and the apparent contradiction completes the argument. 
 \end{rmk}

\begin{rmk}
Theorem \ref{Thm4} provides a somewhat amusing counterpoint to the results in \cite{BBC}. In that work the  non-stickiness condition (\ref{0timeat0}) was used to restore strength and pathwise uniqueness to {\it the degenerate stochastic differential equation 
$$
\mathrm{d} X (t) \,=\, \big| X (t) \big|^\alpha\, \mathrm{d} W (t)
$$
of Girsanov} \cite{G} with $ \, \alpha \in (0, 1/2)\,$ which, in the absence of such a condition, admits several weak solutions. By contrast, Theorem \ref{Thm4} uses the condition  (\ref{0timeat0}) to leave the equation (\ref{sde}) with $\,\lambda>0\,$ bereft of even weak solutions.
\end{rmk}

\begin{rmk}
Consider the equation
\eq
 \label{sdegeq}
\mathrm{d}X(t)\,=\,\kappa \, \mathbf{ 1}_{\{X(t)>0\}}\;\mathrm{d}t + \lambda\, \mathbf{ 1}_{\{X(t)\leq0\}}\;\mathrm{d}t  + \mathbf{ 1}_{\{X(t)\geq0\}}\;\mathrm{d}W(t)\,,\qquad 0 \leq t < \infty 
\en
with diffusion function $\, \sigmab (x) = \mathbf{ 1}_{ [0, \infty)} (x) \,$. In the case $\, \kappa \in \mathbb{R}\,,~\lambda>0\,$ one can follow the lines of the   proof of Theorem \ref{Thm4}, and deduce that there is no weak solution to (\ref{sdegeq}) under the  non-stickiness condition  (\ref{0timeat0}). 

\smallskip
We claim that {\it the equation} (\ref{sdegeq}) {\it fails to have a weak solution also in the case   $\,\kappa = \lambda=0\,$, now even without having to impose the    condition}  (\ref{0timeat0}).  Indeed, by plugging   functions $f\in C_0^\infty(\rr)$ with compact  support  in $(-\infty,0)$, into the Fokker-Planck equation corresponding to the stochastic differential equation \eqref{sdegeq}, we see that  
\[
\pp(X(t)<0)\,=\,\pp(X(0)<0)\,,\qquad 0 \leq t < \infty 
\]
holds for every weak solution $(X,W)$ of the equation \eqref{sdegeq}. Thus, on  the one hand, every weak solution $X$ of the equation \eqref{sdegeq} with $X(0)=x\geq 0$, satisfies $X(t)\geq0$ for Lebesgue almost every $\,t\in [0, \infty)\,$ by Fubini's Theorem. On the other hand, combining the latter conclusion and the equation \eqref{sdegeq} with the  P. L\'evy  characterization of Brownian motion  (e.g., Theorem 3.16, page 157 in \cite{KS}), we conclude that $X$ must be a standard Brownian motion; this  is clearly a contradiction. 

\smallskip
Finally, in the case $\lambda<0$, one can proceed as in section 2 to construct the unique strong solution of the equation 
\eq
 \label{sdegeq2}
\mathrm{d}X(t)\,=\,  \lambda\,  \mathrm{d}t  + \mathbf{ 1}_{\{X(t)\geq0\}}\;\mathrm{d}W(t)\,,\qquad 0 \leq t < \infty\,.
\en  
\end{rmk}


\subsection{Brownian Perturbations that Restore  Strength}

The addition of a suitably correlated Brownian motion   with sufficiently high variance into (\ref{sde}), can restore a pathwise unique, strong solution to this equation  when $\, \lambda >0\,$. Our next result explains how, and its proof works just as well for every value   $\, \lambda \in \mathbb{R}  \,$. 

\begin{thm}
 \label{Thm5}
For any real constant  $\, \lambda  \,$,  and with $\, W$ and $\,V$    standard Brownian motions, the ``perturbed one-sided Tanaka equation"
\eq
\label{sde1}
\mathrm{d}X(t)\,=\,\lambda\;\mathrm{d}t + \mathbf{ 1}_{\{X(t)>0\}}\;\mathrm{d}W(t)+ ( \eta / 2) \,\mathrm{d}V(t) \,,\qquad  0 \leq t < \infty \,,
\en
  has a pathwise unique  strong solution, provided     
     either
  
   \smallskip
    \noindent (i) $\, ~\eta \notin [-1, 1]\,$ and $\, \langle W, V \rangle (t) = - (t / \eta)\,$, $ \, 0 \leq t < \infty\,$, ~or 
  
     \smallskip
     \noindent (ii) $\, \eta \neq 0\,$ and $\, W$, $V$ are independent. 
  
\end{thm}

\noindent{\it Proof:} It is fairly straightforward that solving (\ref{sde1}) under the stated conditions amounts to solving the so-called ``perturbed Tanaka       equation" 
\eq
\label{Prokaj}
\mathrm{d}X(t)\,=\,\lambda\;\mathrm{d}t + \mathrm{sgn}\big( X(t)\big)\;\mathrm{d}M(t)+  \mathrm{d}N(t) \,,\qquad 0 \leq t < \infty \,,
\en
where the processes $\, M  := W / 2\,$, $\,N :=  \big( W  + \eta \, V \big)/ 2\,$ are continuous, orthogonal martingales with quadratic variations $\, \langle M \rangle (t) = t / 4  \,$ and $\, \langle N \rangle (t) = ( \eta^2 - 1) \,t / 4  \,$, respectively. Thus, by the  P. L\'evy theorem once again,  these are independent Brownian motions with respective variance parameters 1/4 and $\,  ( \eta^2 - 1)   / 4  \,$.

The recent work of Prokaj \cite{P} shows that pathwise uniqueness holds for the equation (\ref{Prokaj}). Thus, to complete the proof, it is enough to show that (\ref{Prokaj}) admits a weak solution; for then the Yamada-Watanabe theory (e.g., Corollary 3.23, Chapter 5 in \cite{KS}) guarantees that this solution is actually strong, that is, for all $\, t \in [0, \infty)\,$ we have $$\, \mathcal{F}^X_t \subseteq \mathcal{F}^{\,(M, N)}_t \equiv  \mathcal{F}^{\,(W, V)}_t\,.$$   

In order to prove existence of a weak solution for (\ref{Prokaj}), it is enough to consider the case $\, \lambda =0\,$; this is because a Cameron-Martin-Girsanov change of measure takes then care of any $\, \lambda \in  \mathbb{R}\,$. Therefore, all we need to  do is consider two independent Brownian motions $\,U\,$ and $\,N\,$ with  variance parameters 1/4 and $\,  ( \eta^2 - 1)   / 4  \,$, respectively, along with a  real-valued random variable $\, \zeta\,$ independent of the vector $\, (U, N)$, and define 
\[
X(t)\, :=\, \zeta + U(t) + N(t)\,, \qquad M(t) \,:=\, \int_0^t \mathrm{sgn} \big( X(s) \big)\, \mathrm{d} U(s)\,, \qquad 0 \leq t < \infty\,.
\]
The process $\,M\,$   is a continuous martingale  that satisfies  $$\, \langle M, N \rangle (t ) = \int_0^t \mathrm{sgn}  ( X(s)  ) \, \mathrm{d}  \langle U, N \rangle (s)=0\, \qquad  \mathrm{and}  \qquad \, \langle M  \rangle (t ) =   \langle U \rangle (t ) = t / 4\,;$$ thus, by the P. L\'evy characterization  once again, $\,M\,$ is Brownian motion with variance parameter 1/4, and is independent of the Brownian motion $\,N\,$. But then we have also $\, U(t)  = \int_0^t \mathrm{sgn}  ( X(s) )  \, \mathrm{d} M(s)\,,~~t \ge 0\,$, therefore the representation 
 \[
X(t)  \,=\, \zeta + \int_0^t \mathrm{sgn}  ( X(s) ) \, \mathrm{d} M(s) + N(t)\,,   \qquad 0 \leq t < \infty
\]
as in (\ref{Prokaj}) with $\, \lambda =0\,$; this completes the proof   under the conditions {\it (i)}.

  \smallskip
  Under the conditions of {\it (ii),} the pathwise uniqueness of (\ref{sde1}) is a consequence of Theorem 8.1 in Fernholz et al. \cite{FIKP}, whereas weak existence follows from the results of Bass \& Pardoux \cite{BP}.  \ep
 
 \medskip

\section{Acknowledgments}

We thank Professor Vilmos Prokaj and Dr.$\,\,$Johannes Ruf for reading an early version of the manuscript and offering valuable comments. Extensive discussions with Professor  Hans-J\"urgen Engelbert  and with Drs.$\,$Vasileios Papathanakos and Jon Warren helped sharpen our thinking and arguments.

This research was partially supported by the National Science Foundation under  grants NSF-DMS-08-06211 and  NSF-DMS-09-05754.

 \medskip

\bibliographystyle{alpha}

\begin{thebibliography}{50}

\bibitem{ABTWW} 
\textsc{Appuhamillage, T., Bokil, V., Thomann, E., Waymire, E.} \& \textsc{Wood, B.} (2011) Occupation and local times for skew Brownian motion with application to dispersion along an interface. 
\textit{Annals of Applied Probability} {\bf 21}, 183-214. 

\bibitem{B} 
  \textsc{Barlow, M.T.} (1988)  
 Skew Brownian motion and a one-dimensional stochastic differential equation.
  {\it Stochastics} {\bf 25}, 1-2.
  
  \bibitem{BBC} 
  \noindent
\textsc{Bass, R.F., Burdzy, K.} \& \textsc{Chen, Z.Q.} (2007) Pathwise uniqueness for a degenerate stochastic differential equation. \textit{Annals of Probability} {\bf 35}, 2385-2418. 
  
   \bibitem{BP} 
    \noindent     \textsc {Bass, R. \& Pardoux, E.}  (1987)   Uniqueness of diffusions with piecewise constant co\"efficients.    {\it  Probability Theory and Related Fields}  {\bf 76}, 557-572. 
  
\bibitem{chi}
\textsc{Chitashvili, R. J.} (1997)   
On the nonexistence of a strong solution in the boundary problem for a sticky Brownian motion. {\it Proc. A. Razmadze Math. Inst.} {\bf 115} 17-31.  (Available in preprint form as {\sl CWI Technical Report BS-R8901} (1989),  Centre for Mathematics and Computer Science, Amsterdam.)
  



\bibitem{ES} 
\textsc{Engelbert, H.J.} \& \textsc{Schmidt, W.} (1985) On solutions of stochastic differential equations with no drift. \textit{Zeitschrift f\"ur Wahrscheinlichkeitstheorie und verwandte Gebiete} {\bf 68}, 287-317. 

  \bibitem{FIKP} 
     \textsc{Fernholz, E.R., Ichiba, T., Karatzas, I. \& Prokaj, V.}  (2011) 
   Planar diffusions with rank-based characteristics, and  perturbed Tanaka equations.    {\it Preprint, \textsc{Intech} Investment Management, Princeton.} 

\bibitem{GS} 
\textsc{Gihman, I.I.} \& \textsc{Skorohod, A.V.} (1972)
\textit{Stochastic Differential Equations.}   
Springer Verlag, New York.

\bibitem{G} 
  \textsc{Girsanov, I.V.} (1962)  
 An example of non-uniqueness of the solution to the  stochastic differential equation of K.$\,$It\^o. {\it Theory of Probability and Its Applications} {\bf 7}, 325-331.
 
\bibitem{H}
\textsc{Hajek, B.} (1985) Mean stochastic comparison of diffusions. \textit{Zeitschrift f\"ur Wahrscheinlichkeitstheorie und verwandte Gebiete} \textbf{68}, 315-329.

\bibitem{HS} 
\textsc{Harrison, J.M.} \& \textsc{Shepp, L.A.} (1981) On skew Brownian motion. \textit{Annals of Probability} {\bf 9}, 309-313. 
 


\bibitem{KS} 
\textsc{Karatzas, I.} \& \textsc{Shreve, S.E.} (1991)
\textit{Brownian Motion and Stochastic Calculus.} Second Edition, 
Springer Verlag, New York.

\bibitem{L} 
\textsc{Lejay, A.}  (2006) 
  On the constructions of the skew Brownian motion. 
  {\it Probability Surveys} {\bf 3}, 413-466.
  
\bibitem{MS} 
    \textsc{Manabe, Sh. \& Shiga, T.}  (1973) 
  On   one-dimensional stochastic   differential equations with non-sticky boundary conditions.    {\it  Journal of Mathematics Kyoto University} {\bf 13},  595-603.
  
  \bibitem{N} 
    \textsc{Nakao, S.}  (1972) 
  On the pathwise uniqueness of solutions of one-dimensional stochastic   differential equations.    {\it Osaka Journal of Mathematics} {\bf 9},  513-518.

\bibitem{P} 
   \textsc{Prokaj, V.}  (2010) 
  The solution  of the perturbed Tanaka equation is pathwise unique.  
  {\it Preprint, E\"otv\"os Lor\'and University, Budapest.} 
  
  \bibitem{RY} 
\textsc{Revuz, D.} \& \textsc{Yor, M.} (1999)
\textit{Continuous Martingales and Brownian Motion.} Third Edition, 
Springer Verlag, New York.

 \bibitem{W}
 \textsc{Walsh, J.B.} (1978) 
  A diffusion with   discontinuous local time. 
  In {\it Temps Locaux. Ast\'erisque} {\bf 52-53},  37-45. 
  
   \bibitem{WJ1}
 \textsc{Warren, J.} (1997) 
  Branching processes, the Ray-Knight theorem, and sticky Brownian motion. {\it Lecture Notes in Mathematics} {\bf 1655},  1-15. 
  
     \bibitem{WJ2}
 \textsc{Warren, J.} (1999) 
  On the joining of sticky Brownian motion. {\it Lecture Notes in Mathematics} {\bf 1709},  257-266. 
  
  \bibitem{Z} 
\textsc{Zvonkin, A.K.} (1974) A transformation on the state-space of a diffusion process that removes the drift. \textit{Mathematics of   the USSR (Sbornik)} {\bf 22}, 129-149.


\end{thebibliography}

\bigskip

\end{document}